\documentclass[12pt]{amsart}
\usepackage{amssymb}
\usepackage{amsbsy}
\usepackage{amscd}
\advance\textheight by \topskip
\makeatletter
\renewcommand*\subjclass[2][2000]{%
  \def\@subjclass{#2}%
  \@ifundefined{subjclassname@#1}{%
    \ClassWarning{\@classname}{Unknown edition (#1) of Mathematics
      Subject Classification; using '2000'.}%
  }{%
    \@xp\let\@xp\subjclassname\csname subjclassname@#1\endcsname
  }%
}
\def\cal{\mathcal}
\def\Bbb{\mathbb}
\def\frak{\mathfrak}

\newenvironment{pf*}[1]{\proof[#1]}{\endproof}
\newcommand{\rom}{\textup}
\renewcommand{\thesubsection}{\thesection(\@roman\c@subsection)}
\makeatother
\newtheorem{Theorem}[equation]{Theorem}
\newtheorem{Corollary}[equation]{Corollary}
\newtheorem{Lemma}[equation]{Lemma}
\newtheorem{Proposition}[equation]{Proposition}

\newtheorem{Problem}[equation]{Problem}

\theoremstyle{definition}
\newtheorem{Definition}[equation]{Definition}

\makeatletter
\renewcommand\section{\@startsection{section}{1}%
  {\z@}{.7\linespacing\@plus\linespacing}{.5\linespacing}%
  {\reset@font\normalfont\bfseries\centering}}
\makeatother

\theoremstyle{remark}
\newtheorem{Remark}[equation]{Remark}

\newtheorem*{Acknowledgments}{Acknowledgments}
\numberwithin{equation}{section}
\numberwithin{figure}{section}

\newcommand{\thmref}[1]{Theorem~\ref{#1}}
\newcommand{\secref}[1]{\S\ref{#1}}
\newcommand{\lemref}[1]{Lemma~\ref{#1}}
\newcommand{\propref}[1]{Proposition~\ref{#1}}
\newcommand{\remref}[1]{Remark~\ref{#1}}
\newcommand{\defref}[1]{Definition~\ref{#1}}
\newcommand{\corref}[1]{Corollary~\ref{#1}}
\newcommand{\subsecref}[1]{\S\ref{#1}}

\newcommand{\tabref}[1]{Table~\ref{#1}}

\newcommand{\Romnum}[1]{\expandafter\uppercase\expandafter{\romannumeral #1}} 
\newcommand{\C}{{\Bbb C}}
\newcommand{\Z}{{\Bbb Z}}

\newcommand{\R}{{\Bbb R}}

\newcommand{\CP}{\operatorname{\C P}}

\newcommand{\Spin}{\operatorname{\rm Spin}}
\newcommand{\Spinc}{\Spin^{c}}

\newcommand{\spin}{\operatorname{\frak{spin}}}

\newcommand{\rank}{\operatorname{rank}}
\newcommand{\Sign}{\operatorname{Sign}}

\newcommand{\ind}{\mathop{\text{\rm ind}}\nolimits}
\newcommand{\tr}{\operatorname{tr}}

\newcommand{\CC}{{\cal C}} 
\newcommand{\D}{{\cal D}} %
\newcommand{\UU}{{\frak U}} %
\newcommand{\SW}{\operatorname{SW}}

\begin{document}
\title[Nonsmoothable group actions on elliptic surfaces]
{Nonsmoothable group actions on elliptic surfaces}
\author{Ximin Liu}
\address{Department of Applied Mathematics, Dalian University of Technology, Dalian 116024, China}
\email{ximinliu@dl.cn}

\author{Nobuhiro Nakamura}
\address{Graduate School of Mathematical Sciences, University of Tokyo, 3-8-1, Komaba, Meguro-ku, Tokyo, 153-8914, Japan}
\email{nobuhiro@ms.u-tokyo.ac.jp}
\thanks{The first author is supported by NSFC(10771023), RFDP(20050141011) and NCET}
\begin{abstract}
Let $G$ be a cyclic group of order $3$, $5$ or $7$, and $X=E(n)$ be the relatively minimal elliptic surface with rational base. 
In this paper, we prove that under certain conditions on $n$, 
there exists a locally linear $G$-action on $X$ which is nonsmoothable with respect to infinitely many smooth structures on $X$.  This extends the main result of \cite{LN}.
\end{abstract}
\maketitle
%
%
\section{Introduction}\label{sec:intro}
%
%
It is a classical result that every finite group action on a surface is equivalent to a smooth one. 
In higher dimensions, there exist examples of nonsmoothable actions.
Since bad local behavior is often the reason why these actions can not be smooth, one can naturally ask whether {\it locally linear} actions are smoothable or not.
In \cite{KLee}, S.~Kwasik and K.~B.~Lee proved that in dimension 3 a finite group action is smoothable if and only if it is locally linear. 
However, in dimensions higher than $3$, this is not true. 
In fact, many examples of nonsmoothable locally linear actions are known \cite{KLee,KLaw,HL,Bryan,Kiyono,ChenKw}.

The authors also constructed such a nonsmoothable action.
Let $\Z_p$ be the cyclic group of order $p$.
\begin{Theorem}[\cite{LN}]\label{thm:stdk3}
There exists a locally linear pseudofree $\Z_3$-action on a $K3$ surface $X$ which is nonsmoothable with respect to the standard smooth structure on $X$.  
\end{Theorem}
\begin{Remark}
An action on a space is called {\it pseudofree} if it is free on the complement of a finite subset. 
Note that every smooth action is locally linear.
On the other hand  a locally linear action is not necessarily smooth.
\end{Remark}

In this paper, we extend \thmref{thm:stdk3} to two directions.
One direction is to consider elliptic surfaces with higher Euler numbers. 
Another direction is to consider higher order cyclic group actions.
In fact, we will prove the following.
\begin{Theorem}\label{thm:elln}
Let $G=\Z_p$, where $p=3$, $5$ or $7$, and $X=E(n)$ be the relatively minimal simply-connected elliptic surface with rational base which has the Euler number $12n$. 
Suppose $n$ is even and $n\geq 2$, and let
\begin{equation}\label{eq:cn}
c_{n-2} :=
\left(
\begin{gathered}
n-2\\
\frac{n-2}2
\end{gathered}
\right). 
\end{equation}
{\rom(}Assume $c_0=1$.{\rom)}
If  $c_{n-2}\not\equiv 0$ mod $p$, then there exists a locally linear $G$-action on $X$ which is nonsmoothable with respect to infinitely many smooth structures on $X$.
\end{Theorem}
\begin{Remark}
We do not know whether there exists a smooth structure on $X$ on which the above locally linear action is smoothable, or not.
\end{Remark}
\begin{Remark}\label{rem:p-adic}
The number $c_{n-2}$ is the Seiberg-Witten invariant of the standard $E(n)$ for the $\Spinc$-structure $c_{\spin}$ associated to the spin structure.
For a prime $p$, it is  easy to see that the condition $c_{n-2} \not\equiv 0$ mod $p$ is equivalent to the condition that the $p$-adic expansion of $\frac{n-2}2$ does not contain any number bigger than $p/2$ in any tab.
\end{Remark}

The proof of \thmref{thm:elln} is analogous to that of \thmref{thm:stdk3}(\cite{LN}), and it is divided into two steps: 
In the first step, we give a constraint on smooth actions.
In the second step, we construct a locally linear action which would violate the constraint if it were smooth. 

To obtain a constraint on smooth actions, we use the Seiberg-Witten gauge theory.
In fact, we use a mod $p$ vanishing theorem of Seiberg-Witten invariants, which is originally proved by Fang \cite{Fang}, and generalized by the second author \cite{Nakamura}, with known calculations of the Seiberg-Witten invariants of elliptic surfaces.
On the other hand, to construct locally linear actions, we invoke a remarkable realization theorem by Edmonds and Ewing\cite{EE}.

\subsection{$\Z_3$-actions on elliptic surfaces}

First we will explain the case of $\Z_3$-actions more precisely.
When we fix a generator $g$ of $G=\Z_3$, the representation at a fixed point can be described by a pair of nonzero integers $(a,b)$ modulo $3$ which is well-defined up to order and changing the sign of both together.
Hence, there are two types of fixed points.
\begin{itemize}
\item The type ($+$): $(1,2) = (2,1)$.
\item The type ($-$): $(1,1) = (2,2)$.
\end{itemize}
Let $m_+$ be the number of fixed points of the type ($+$), and $m_-$ be the number of fixed points of the type ($-$).

The Euler number of a $4$-manifold $X$ is denoted by $\chi(X)$, and the signature by $\Sign(X)$.
For any $G$-space $V$, let $V^G$ be the fixed point set of the $G$-action. Let $b_{\bullet}^G = \dim H_{\bullet}(X;\R)^G$, where $\bullet = 2,+,-$.

We will prove the following constraint on smooth $G$-actions.
\begin{Theorem}\label{thm:spin-z3}
Let $G=\Z_3$, and $X$ be a simply-connected closed oriented smooth spin $4$-manifold with $b_+\geq 2$, which satisfies $2\chi(X)+3\Sign(X)=0$.
Suppose $G$ acts on $X$ smoothly and pseudofreely so that $b_+^G\geq 1$. 
If the Seiberg-Witten invariant $\SW_X(c_{\spin})$ for the $\Spinc$-structure $c_{\spin}$ associated to the spin structure is not divisible by $3$, then
\begin{equation}\label{eq:m0}
m_+=0\text{ or }m_-=0.
\end{equation}
\end{Theorem}
Note that the spin manifold in \thmref{thm:spin-z3} is a homotopy $E(n)$. 
The Seiberg-Witten invariants of elliptic surfaces have been already calculated \cite{FM2, FS0}.
Later, we will introduce an infinite family of smooth structures on $E(n)$, denoted by $\UU_{E(n),p}$, which has the property that, for $X\in\UU_{E(n),p}$, $\SW_X(c_{\spin})\not\equiv 0$ mod $p$ if $c_{n-2}\not\equiv 0$ mod $p$.
(See \defref{def:up}.)
Such smooth structures are obtained by log transformations and Fintushel-Stern's knot surgery constructions \cite{FS1}.  
For smooth actions on smooth structures in $\UU_{E(n),3}$, the following holds by \thmref{thm:spin-z3}. 
\begin{Corollary}\label{cor:ell-z3}
Let $X$ be a homotopy $E(n)$ of even and positive $n$ with a smooth structure in $\UU_{E(n),3}$. 
Suppose that $G=\Z_3$ acts on $X$ smoothly and pseudofreely so that $b_+^G\geq 1$. 
If $c_{n-2}$ in \eqref{eq:cn} satisfies $c_{n-2}\not\equiv 0$ mod $3$, 
 then $m_+=0$ or $m_-=0$.
\end{Corollary}
\begin{Remark}
If a pseudofree $G$-action on a $K3$ surface $X$ is {\it holomorphic}, then $G$ acts on the space of holomorphic $2$-forms which consists of constant sections of the canonical line bundle. 
Then, it follows that the weight of the $G$-action on the fiber over each fixed point is always same, i.e., $m_+=0$ or $m_-=0$.
\corref{cor:ell-z3} might be considered as a generalization of this fact.
\end{Remark}
On the other hand, by using the result by Edmonds and Ewing \cite{EE}, we can construct locally linear actions on simply-connected manifolds realizing given fixed point data. 
In fact, a locally linear $G$-action with $m_+>0$ and $m_->0$ can be constructed on each $E(n)$.  
The $\Z_3$-case of \thmref{thm:elln} will be proved by summing up these.

In particular, we will give a detailed proof of the following result on $\Z_3$-actions on the homotopy $E(4)$ as a model case.
\begin{Theorem}\label{thm:ellz3}
Let $G=\Z_3$, and $X=E(4)$.
For locally linear pseudofree $G$-actions on $X$, we have the following{\rom :}
\begin{enumerate}
\item  Every locally linear pseudofree $G$-action on $X$ belongs to one of ten classes in \tabref{tab:ellz3}.
Furthermore, each of classes except the class $C_1$ can be actually realized by a locally linear pseudofree $G$-action on $X$.
\begin{table}[h]
\caption{$\Z_3$-actions on $E(4)$}
\label{tab:ellz3}
\begin{center}
\begin{tabular}{l|c|c|c|c|c|c|c|c}
Class & $\#X^G$ & $m_+$ & $m_-$ & $b_2^G$ & $b_+^G$ & $b_-^G$ & $\Sign(X/G)$ &\\
\hline
$A_1$ & $12$ & $12$ & $0$ & $22$ & $7$ & $15$ & $-8$  &\\
$A_2$ & $15$ & $9$ & $6$ & $24$ & $7$ & $17$ & $-10$ & {\footnotesize {\rm NS}}\\
$A_3$ & $18$ &  $6$ & $12$ & $26$ & $7$ & $19$ & $-12$ & {\footnotesize {\rm NS}}\\
$A_4$ & $21$ & $3$ & $18$ & $28$ & $7$ & $21$ & $-14$ &{\footnotesize {\rm NS}} \\
$A_5$ & $24$ &  $0$ & $24$ & $30$ & $7$ & $23$ & $-16$ &\\
\hline
$B_1$ & $9$ & $6$ & $3$ & $20$ & $5$ & $15$ & $-10$ & {\footnotesize {\rm NS}}\\
$B_2$ & $12$ & $3$ & $9$ & $22$ & $5$ & $17$ & $-12$ & {\footnotesize {\rm NS}}\\
$B_3$ & $15$ & $0$ & $15$ & $24$ & $5$ & $19$ & $-14$ &\\
\hline
$C_1$ & $3$ & $3$ & $0$ & $16$ & $3$ & $13$ & $-10$ &\\
$C_2$ & $6$ & $0$ & $6$ & $18$ & $3$ & $15$ & $-12$ &\\
\end{tabular}
\end{center}
\end{table}
\item Every locally linear action in the classes $A_2$, $A_3$, $A_4$, $B_1$ and $B_2$ {\rom (}indicated as ``NS'' in \tabref{tab:ellz3}{\rom )} is nonsmoothable for infinitely many smooth structures in $\UU_{E(4),3}$. 
\end{enumerate}
\end{Theorem}
\begin{Remark}
In \thmref{thm:ellz3} (and also in \thmref{thm:k3z5}), 
the classification into classes enumerates all candidates of fixed point data which satisfy the Lefschetz formula and the $G$-signature formula.  
\end{Remark}
\begin{Remark}
At present, we do not know any concrete example of smooth $\Z_3$-actions on $E(4)$ for any smooth structure.
%
%
In particular,
we do not know whether there exists a smooth structure on which a class indicated as ``NS'' ($A_2$, $A_3$, $A_4$, $B_1$ or $B_2$) can be realized by a smooth action, or not.
\end{Remark}
\begin{Remark}
We can obtain similar classification results on $\Z_3$-actions on homotopy $E(n)$ for larger $n$. 
The classifications for $E(8)$ and $E(10)$ are summarized in \tabref{table:e8-z3} and \tabref{table:e10-z3} at the end of the paper.
See also \remref{rem:en-z3}.
\end{Remark}
%
%
\subsection{$\Z_5$-actions on $K3$}
%
%
For $\Z_5$-actions on $K3$, the following holds.

\begin{Theorem}\label{thm:k3z5}
Let $G=\Z_5$.
For locally linear pseudofree $G$-actions on a $K3$ surface $X$, the following hold{\rom :}
\begin{enumerate}
\item  Locally linear pseudofree $G$-actions on $X$ are classified into $285$ classes. 
Furthermore, each of $285$ classes can be actually realized by a locally linear pseudofree $G$-action on $X$.
\item There are $41$ classes which can not be realized by a smooth action with respect to infinitely many smooth structures in $\UU_{E(2),5}$. 
\end{enumerate}
\end{Theorem}
\begin{Remark}
At least, two classes in the total $285$ classes can be realized by smooth actions.
See \propref{prop:z5-1} and \remref{rem:z5-502}.
The authors do not know whether the other classes can be realized by smooth actions, or not.
\end{Remark}
Before concluding the introduction, we would like to mention a result by W.~Chen and S.~Kwasik.
In the paper \cite{ChenKw}, Chen and Kwasik proved the existence of a family of symplectic exotic $K3$'s on which every nontrivial ${\Bbb Z}_p$-action of prime $p\geq 7$ is nonsmoothable.
It would be interesting to compare the method of Chen-Kwasik with our approach.
Chen-Kwasik's argument uses
\begin{itemize}
\item the fact that the Seiberg-Witten basic classes are preserved by symmetries, and 
\item techniques in symplectic geometry to investigate symmetries of symplectic manifolds.
\end{itemize}
In fact, a key point is that the exotic $K3$ they constructed have many basic classes, and this fact with symplectic techniques implies a strong restriction on smooth actions on them. 
Note that their method can not be applied to actions on the standard $K3$ whose only Seiberg-Witten basic class is $0$. 
Moreover, although there is some possibility of applying Chen-Kwasik's method to some kind of exotic elliptic surfaces, we do not succeed in constructing a nonsmoothable action on a standard elliptic surface by their method, because a standard elliptic surface does not have sufficiently many basic classes.
Thus, the method of Chen-Kwasik is ``perpendicular'' to ours. 
(See also \subsecref{subsec:dep}.)

The paper is organized as follows. 
Section 2 provides preliminaries.
Sections 3, 4 and 5 deal with $\Z_3$, $\Z_5$ and $\Z_7$-actions on elliptic surfaces, respectively. 
In Section 6, some concluding remarks are given.
\begin{Acknowledgments}
The authors would like to thank M.~Furuta for invaluable discussions. 
The second author would like to thank K.~Oguiso for teaching him on log Enriques surfaces and recommending the paper \cite{Zh}. 
It is also a pleasure to thank R.~Stern for valuable comments, W.~Chen and S.~Kwasik for sending their preprints, and the referee for valuable suggestions. 
\end{Acknowledgments}
%
%
%
\section{Preliminaries}\label{sec:pre}
%
%
The purpose of this section is to collect known facts on smooth and locally linear actions.
%
%
%
\subsection{The $G$-index theorems}\label{subsec:index}
%
%
For the generality of $G$-index theorems, we refer \cite{AB,AB2, AH, Sh}.
Let $G=\Z_p$ of prime $p$, and fix a generator $g$. 
Suppose $G$ acts on a closed smooth $4$-manifold $X$ smoothly and {\it pseudofreely}, and the fixed point data for the generator $g$ are given as $\{(a_i,b_i)\}_{i=1}^{N}$.

The {\it $G$-signature formula }is,
\begin{equation}\label{eq:G-sgn0}
\Sign(g, X) = \sum_{i=1}^{N}s_{a_ib_i}
\end{equation}
where 
\begin{equation}\label{eq:sxy}
s_{xy}=\frac{(\zeta^{x} + 1)(\zeta^{y} + 1)}{(\zeta^{x} - 1)(\zeta^{y} - 1)},
\end{equation}
and $\zeta = \exp(2\pi\sqrt{-1}/p)$.

Suppose further that $X$ is {\it spin} and the $G$-action is a {\it spin} action. 
Let $D_X$ be the $G$-equivariant Dirac operator.
Then the {\it $G$-spin theorem} is,
\begin{equation}\label{eq:G-spin}
\ind_g D_{X} = \sum_{i=1}^N p_{a_ib_i},
\end{equation}
where 
\begin{equation}\label{eq:pxy}
p_{xy} = \frac1{{(\zeta^{x})}^{1/2} - {(\zeta^{x})}^{-1/2}}\frac1{{(\zeta^{y})}^{1/2} - {(\zeta^{y})}^{-1/2}},
\end{equation}
and signs of ${(\zeta^{x})}^{1/2}$ and ${(\zeta^{y})}^{1/2}$ are determined by the rule
$$
\left\{{(\zeta^{x})}^{1/2}\right\}^p =\left\{{(\zeta^{y})}^{1/2}\right\}^p = 1.  
$$
(This is because $G$ is supposed odd order, and the $g$-action on the spin structure generates the $G$-action on the spin structure. See \cite[p.20]{AH} or \cite[p.175]{Sh}.)
%
\subsection{The realization theorem by Edmonds and Ewing}\label{subsec:EE}
We summarize the realization theorem of locally linear pseudofree actions by Edmonds and Ewing \cite{EE} in the special case when $G=\Z_3, \Z_5$ or $\Z_7$.
\begin{Theorem}[\cite{EE}]\label{thm:EE}
Let $G=\Z_p$, where $p=3,5$ or $7$.
Suppose that one is given a fixed point data
$$
\D = \{(a_0,b_0), (a_1,b_1), \ldots, (a_n,b_n), (a_{n+1},b_{n+1})\},
$$
where $a_i, b_i \in \Z_p\setminus\{0\}$, and a $G$-invariant symmetric unimodular form
$$
\Phi\colon V\times V\to \Z,
$$
where $V$ is a finitely generated $\Z$-free $\Z[G]$-module.
Then the data $\D$ and the form $(V,\Phi)$ are realizable by a locally linear, pseudofree, $G$-action on a closed, simply-connected, topological $4$-manifold if and only if they satisfy the following two conditions{\rom :}
\begin{enumerate}
\item The condition REP{\rom :} As a $\Z[G]$-module, $V$ splits into $F\oplus T$, where $F$ is free and $T$ is a trivial $\Z[G]$-module with $\rank_\Z T = n$.
\item The condition GSF{\rom :} The $G$-Signature Formula is satisfied{\rom :}
\begin{equation}\label{eq:G-sgn}
\Sign(g, (V,\Phi)) = \sum_{i=0}^{n+1}\frac{(\zeta^{a_i} + 1)(\zeta^{b_i} + 1)}{(\zeta^{a_i} - 1)(\zeta^{b_i} - 1)},
\end{equation}
where $\zeta = \exp(2\pi\sqrt{-1}/p)$.
\end{enumerate}
\end{Theorem}
\begin{Remark}
In \cite{EE}, Edmonds and Ewing proved the realization theorem for all cyclic groups of prime order $p$. 
For general  $p$, the third condition {\it TOR} which is related to the Reidemeister torsion should be satisfied.
However, when $p$ is a prime less than $23$, the condition {\it TOR} is redundant.
This follows from the fact that the class number of $\Z[\zeta]$ is $1$, and Corollary 3.2 of \cite{EE}.
\end{Remark}
%
%
\subsection{Mod $p$ vanishing theorem of Seiberg-Witten invariants}
%
Let $p$ be an odd prime, and suppose that $G=\Z_p$ acts {\it smoothly} on a smooth closed oriented $4$-manifold $X$ with $b_1=0$, $b_+\geq 2$. 
Fix a $G$-invariant metric.
Suppose that the $G$-action lifts to a $\Spinc$-structure $c$.
Fix a $G$-invariant connection $A_0$ on the determinant line bundle $L$ of $c$. Then the Dirac operator $D_{A_0}$ associated to $A_0$ 
is $G$-equivariant, and 
the $G$-index of $D_{A_0}$ can be written as $\ind_{G} D_{A_0} = \sum_{j=0}^{p-1} k_j\C_j\in R(G)\cong\Z [t]/(t^p-1)$, where $\C_j$ is the complex $1$-dimensional weight $j$ representation of $G$ and $R(G)$ is the representation ring of $G$.

In such a situation, the following theorem is proved.
\begin{Theorem}[\cite{Fang,Nakamura}]\label{thm:mod-p}
Suppose further that $b_+^G:=\dim H^+(X;\R)\geq 1$.
If $2k_j < 1 + b_+^G$ for $j=0,1,\ldots ,p-1$, then the Seiberg-Witten invariant $\SW_X(c)$ for $c$ satisfies
$$
\SW_X(c)\equiv 0 \mod p.
$$
\end{Theorem}
\begin{Remark}
In \cite{Fang}, Fang suppose that $b_+^G=b_+$. 
In \cite{Nakamura}, the second author weakened that condition as above, and generalized to the case when $b_1\geq 1$.
\end{Remark}
\begin{Remark}\label{rmk:spin-action}
Suppose $X$ is spin and simply-connected.
Let $c_{\spin}$ be the $\Spinc$-structure associated to the spin structure, whose determinant line bundle $L$ is trivial.
If $p$ is odd, then every $G=\Z_p$-action on $X$ has a spin lift. 
Therefore it has a lift to $c_{\spin}$ such that there exists a trivialization $L=X\times\C$ of the determinant line bundle, and the induced $G$-action on $L$ is given by the diagonal action of the $G$-action on $X$ and the trivial action on $\C$.
Let $A_0$ be the trivial flat $G$-invariant connection on $L$. 
Then $\ind_g D_{A_0}$ can be calculated by the $G$-spin theorem \eqref{eq:G-spin}.
\end{Remark}
%
%
%
%
%
\section{$\Z_3$-actions on elliptic surfaces}\label{sec:z3}
%
%
%
In this section, we prove the $\Z_3$-case of our main theorem(\thmref{thm:elln}).
Although the final goal is to give a proof for $E(n)$ of general $n$, we first prove \thmref{thm:ellz3} on $E(4)$ as a model case.
In the course of the proof, \thmref{thm:spin-z3} and \corref{cor:ell-z3} are also proved.
%
%
\subsection{Existence of locally linear $\Z_3$-actions}\label{subsec:loclin-z3}
%
%
Let us begin the proof of the assertion (1) of \thmref{thm:ellz3}.
Let $X$ be $E(4)$ in \thmref{thm:ellz3}.
Suppose that a locally linear pseudofree $G$-action on $X$ is given.
Let $e=\chi(X)$ and $s=\Sign(X)$.
First of all, the ordinary Lefschetz formula should hold: $L(g,X) = 2 + \tr ( g|_{H^2(X)}) =\#X^G$.
Since $\#X^G = m_+ + m_- $ and $ 2 + \tr ( g|_{H^2(X)}) \leq e$, we obtain
\begin{equation}\label{eq:ineq-e}
m_+ + m_- \leq e.
\end{equation}
%
%
Note that
\begin{equation*}
\chi (X/G) = \frac13 \{ e + 2 (m_+ + m_- )\}.
\end{equation*}
Since $\chi (X/G)$ is an integer, we have
\begin{equation}\label{eq:mod3}
m_+ + m_- \equiv -\frac12 e \mod 3.
\end{equation}

By \thmref{thm:EE}, the $G$-Signature Formula should hold:
\begin{align*}
\Sign (g,X) &= \Sign (g^2,X) = \frac13 (m_+ - m_-),\\
\Sign (X/G) &=  \frac13 \left\{ s + \frac23 (m_+ - m_-)\right\}.\\
\end{align*}
Since $\Sign (X/G)$ is an integer, 
\begin{equation}\label{eq:mod9-s}
m_+ - m_- \equiv -\frac32 s \mod 9.
\end{equation}

We can calculate $b_+^G$ and $b_-^G$ from $\chi(X/G)$ and $\Sign (X/G)$:
\begin{align}\label{eq:bpG}
b_+^G &= \frac16\left\{e+s+\frac13(8m_++ 4m_-)\right\}-1,\\
b_-^G &= \frac16\left\{e-s+\frac13(4m_++ 8m_-)\right\}-1.
\end{align}
These should satisfy
\begin{equation}\label{eq:ineq-bG}
0\leq b_+^G \leq b_+,\, 0\leq b_-^G \leq b_-.
\end{equation}
By \eqref{eq:ineq-e}, \eqref{eq:mod3}, \eqref{eq:mod9-s}, \eqref{eq:ineq-bG} and non-negativity of $m_+$ and $m_-$, we obtain \tabref{tab:ellz3}.

Note that the above argument proves the following result for more general $X$.
\begin{Proposition}\label{prop:propz3}
Suppose that $G=\Z_3$ acts on a simply-connected closed oriented manifold $X$ locally linearly and pseudofreely.
Let $e=\chi(X)$ and $s=\Sign(X)$. 
Then the data $(m_+, m_-)$ satisfies \eqref{eq:ineq-e}, \eqref{eq:mod3}, \eqref{eq:mod9-s} and \eqref{eq:ineq-bG}.
\end{Proposition}

Next we will prove the existence of actions. 
To prove the existence of locally linear actions, we invoke \thmref{thm:EE}. 
We need to construct $G$-actions on the intersection form.
Let $(V_X, \Phi_X)$ be the intersection form of $X=E(4)$. 
Since an even indefinite form is completely characterized  by  its rank and signature, $(V_X, \Phi_X)$ is isomorphic to $7H\oplus 2\Gamma_{16}$, where $H$ is the hyperbolic form, and $\Gamma_{16}$ is a negative definite even form of rank $16$ given below. 
We will construct $G$-actions on $3H$ and $\Gamma_{16}$ separately.

Let $r$ be a multiple of $4$, and  $\Gamma_{r}$ be the lattice of $(x_1,\ldots,x_{r})\in(\frac12 \Z)^{r}$ which satisfy
\begin{enumerate}
\item $x_i\equiv x_j\mod \Z$ for any $i,j$,
\item $\sum_{i=1}^{r}x_i\equiv 0\mod 2\Z$.
\end{enumerate}
The unimodular bilinear form on $\Gamma_{r}$ is defined by $-\sum_{i=1}^{r}x_i^2$. 
Then $\Gamma_{r}$ is even and negative-definite.

\begin{Lemma}\label{lem:even}
Let $r=16(3q+1)$, where $q$ is a non-negative integer.  
For each integer $k$ which satisfies $0\leq k\leq 16q+5$, there is a $G$-action on $\Gamma_{r}$ such that
$$
\Gamma_{r}\cong (r-3k)\Z\oplus k\Z[G]\text{ as a $\Z[G]$-module}.
$$
This unimodular form with the $G$-action is denoted by $\Gamma_{r,k}$.
\end{Lemma}
\proof
When $k=0$, it suffices to take the trivial $G$-action. 
Hence we suppose $k\geq 1$. 

Note that the symmetric group of degree $r$ acts on $\Gamma_{r}$ as permutations of components.
For the fixed generator $g$ of $G$, define the $G$-action on $\Gamma_{r}$ by
$$g = (1,2,3)(4,5,6)\cdots(3k-2,3k-1,3k),$$
where $(l,m,n)$ is the cyclic permutation of $(x_l, x_m ,x_n)$.

As a basis for $\Gamma_{r}$, we take 
\begin{equation*}
f_i = \left\{
\begin{aligned}
e_i+e_{r},& \quad\qquad (i=1,\ldots ,24q+9), \\
e_i-e_{r},& \quad\qquad (i=24q+10,\ldots,r-1),\\
\frac12 (e_1+e_2&+\cdots+e_{r}),\, (i=r),
\end{aligned}\right.
\end{equation*}
where $e_1,\ldots,e_{r}$ is the usual orthonormal basis for $\R^{r}$. 
Then the basis $(f_1,f_2,\ldots,f_{r})$ gives the required direct splitting.
\endproof
For a $G$-form $F$, $r_+^G$ (resp. $r_-^G$) denote the rank of the $G$-fixed part of a maximal positive (resp. negative) definite subspace of $F\otimes\R$.
\begin{Lemma}\label{lem:H}
There exist the following $G$-invariant hyperbolic forms.
\begin{enumerate}
\item $A$ such that $A \cong H$ as a form and $A\cong\Z\oplus\Z$ as a $\Z[G]$-module.
\item $B_{2,0}$ such that $B_{2,0} \cong 2H$ as a form and $B_{2,0}\cong\Z\oplus\Z[G]$ as a $\Z[G]$-module and $r_+^G=2$ and $r_-^G=0$.
\item $B_{0,2}$ such that $B_{0,2} \cong 2H$ as a form and $B_{0,2}\cong\Z\oplus\Z[G]$ as a $\Z[G]$-module and $r_+^G=0$ and $r_-^G=2$.
\item $C_{1,1}$ such that $C_{1,1} \cong 3H$ as a form and $C_{1,1}\cong\Z[G]\oplus\Z[G]$ as a $\Z[G]$-module and $r_+^G=r_-^G=1$.
\end{enumerate}
\end{Lemma}
\proof
(1) is trivial.
(4)The form $C_{1,1}$ is given as permutations of three $H$'s.

(2)
Let us consider the basis of  $\Z\oplus\Z[G]$ of the form $\{f,e,ge,g^2e\}$. 
With respect to this basis, let us consider the form represented by the matrix
$$
P=
\begin{pmatrix}
 2 &-1 &-1 &-1\\ 
-1 & 0 & 1 & 1\\
-1 & 1 & 0 & 1\\
-1 & 1 & 1 & 0
\end{pmatrix}.
$$
It is easy to see that $P$ represents a unimodular even form on $\Z\oplus\Z[G]$ and that $r_+^G=2$ and $r_-^G=0$. 
Since the rank of $P$ is $4$, $P$ should be equivalent to $2H$. 

(3)The form represented by $-P$ is the required one.
\endproof

With \lemref{lem:even} and \lemref{lem:H} understood, we can construct required $G$-invariant unimodular forms for classes on $E(4)$ in \thmref{thm:ellz3}:
\begin{itemize}
\item The class $A_1$: $3A\oplus 2B_{2,0}\oplus \Gamma_{16,5}\oplus\Gamma_{16,5}$.
\item The class $A_2$: $5A\oplus B_{2,0}\oplus \Gamma_{16,5}\oplus\Gamma_{16,5}$.
\item The class $A_3$: $7A\oplus \Gamma_{16,5}\oplus\Gamma_{16,5}$.
\item The class $A_4$: $7A\oplus \Gamma_{16,5}\oplus\Gamma_{16,4}$.
\item The class $A_5$: $7A\oplus \Gamma_{16,5}\oplus\Gamma_{16,3}$.
\item The class $B_1$: $A\oplus 2B_{2,0}\oplus B_{0,2} \oplus\Gamma_{16,5}\oplus\Gamma_{16,5}$.
\item The class $B_2$: $3A\oplus B_{2,0}\oplus B_{0,2}\oplus \Gamma_{16,5}\oplus\Gamma_{16,5}$.
\item The class $B_3$: $5A\oplus B_{0,2} \oplus\Gamma_{16,5}\oplus\Gamma_{16,5}$.
\item The class $C_2$: $B_{2,0}\oplus B_{0,2} \oplus C_{1,1}\oplus\Gamma_{16,5}\oplus\Gamma_{16,5}$.
\end{itemize}

By our method, we can not construct a $G$-form for the class $C_1$ in \thmref{thm:ellz3}. 

Now, for each class above, the conditions {\it REP} and {\it GSF} are satisfied. 
Therefore, by \thmref{thm:EE}, there exists a  closed simply-connected $4$-manifold $X^\prime$ with a locally linear pseudofree $G$-action realizing each given data.
Since $X^\prime$ is simply-connected and has even intersection form, we see that $X^\prime$ is homeomorphic to $E(4)$ by Freedman's theorem \cite{Freedman}.
Thus the assertion (1) of \thmref{thm:ellz3} is proved.
%
%
%
\subsection{A constraint on smooth $\Z_3$-actions on elliptic surfaces}\label{subsec:const-z3}
%
%
%
First, we prove \thmref{thm:spin-z3}.

\proof[Proof of \thmref{thm:spin-z3}]
Let $G$ act on $c_{\spin}$ as in \remref{rmk:spin-action}.
Take the trivial flat connection $A_0$ as the reference $G$-invariant connection.
By \thmref{thm:mod-p}, $\SW_X(c_{\spin})\not\equiv 0$ mod $3$ implies that there exist $j$ which satisfies $2k_j\geq 1+b_+^G$.

Note that $b_+^G$ is calculated in \eqref{eq:bpG}. 
Coefficients $k_j$ are calculated by the $G$-spin theorem. 
By the $G$-spin theorem \eqref{eq:G-spin}, we have
\begin{align*}
\ind_g D_{A_0} &= k_0 + \zeta k_1 + \zeta^2 k_2 = \frac13(m_+-m_-),\\
\ind_{g^2} D_{A_0} &= k_0 + \zeta^2 k_1 + \zeta k_2 = \frac13(m_+-m_-),\\
\ind_1 D_{A_0} &= k_0 + k_1 + k_2 = - \frac18 s.
\end{align*}
Solving these, we have
\begin{align*}
k_0 &= \frac29(m_+-m_-) - \frac1{24}s,\\
k_1=k_2&=-\frac19(m_+-m_-) - \frac1{24}s.
\end{align*}
>From these and the relation $2e+3s=0$, we have $m_+=0$ or $m_-=0$. 
\endproof
\begin{Remark}
We have an example of smooth $G$-action on a spin manifold $X$ which satisfies the assumption of \thmref{thm:spin-z3}, however $\SW_X(c_{\spin})\equiv 0$ mod $3$ and $m_+>0$ and $m_->0$. 
Let us consider the lattice $\Z\oplus\zeta\Z\subset \C$, where $\zeta=\exp (2\pi\sqrt{-1}/3)$, and let $T$ be the  $2$-torus $\C/(\Z\oplus\zeta\Z)$ on which $G$ acts by multiplications by $\zeta$. 
Consider a $2$-sphere $S$ with a $G$-action, where the $G$-action is generated by the $2\pi/3$-rotation. 
Let $N$ be $S\times T$ with the diagonal $G$-action. 
We consider the projection $N\to S$ as an elliptic fibration with a $G$-action.
Choosing a free point $q$ on $S$, and taking fiber connected sum of $N$ with $3$ Kummer surfaces over three points $q$, $gq$ and $g^2q$, we obtain $E(6)$ with a $G$-action. 
For this $G$-action on $E(6)$, $m_+=m_-=3>0$ and $\SW_X(c_{\spin})=6$.  
\end{Remark}
Now, we discuss smooth structures on elliptic surfaces and their Seiberg-Witten invariants.
There are two well-known methods to produce exotic smooth structures on $E(n)$: logarithmic transformations and Fintushel-Stern's knot surgery construction.

Logarithmic transformations produce elliptic surfaces with multiple fibers. 
(See e.g. \cite{GS}.)
Let us consider the case of at most $2$ multiple fibers $E(n)_{k,l}$. 
(Here, we assume $k$ and $l$ may be $1$.)
Suppose that $n$ is even and positive, $k$ and $l$ are odd, and $\gcd(k,l)=1$. 
Note that these conditions imply that $E(n)_{k,l}$ is spin and simply-connected.
For such  $E(n)_{k,l}$, the following are known:
\begin{enumerate}
\item $E(n)_{k,l}$ is homeomorphic to $E(n)$ 
if and only if $k$ and $l$ are odd and mutually coprime. (See e.g. \cite{Ue}.)
\item $E(n)_{k,l}$ is diffeomorphic to $E(n)_{k^\prime,l^\prime}$ if and only if $\{k,l\}=\{k^\prime,l^\prime \}$ as unordered pair. (See \cite{GS}.)
\item $E(2)=E(2)_{1,1}$ (no multiple fiber) is diffeomorphic to the standard $K3$ surface.
\end{enumerate}
To prove (2), Seiberg-Witten invariants are used.
In particular, the Seiberg-Witten invariant of $E(n)_{k,l}$ for $c_{\spin}$ is given by 
\begin{equation}\label{eq:ell-spin}
\SW_{E(n)_{k,l}}(c_{\spin}) = (-1)^{\frac{ n-2}2}
\left(
\begin{gathered}
n-2\\
\frac{n-2}2
\end{gathered}
\right).
\end{equation}
Note that this is independent on $k,l$. (See \cite{FM2, FS1}.)

Fintushel and Stern introduced the knot surgery construction in \cite{FS1}, which enables us to produce more exotic smooth structures on $E(n)$.
The construction is given as follows. 
(See \cite{FS1} for details.)
For each $n$, $X=E(n)$ admits a elliptic fibration which contains a cusp fiber. 
Take a smooth embedded torus $T$ in a regular neighborhood of a cusp fiber which represents a nontrivial homology class.
Remove a tubular neighborhood of $T$ from $X=E(n)$, and denote the resulting manifold by $X^\prime$. 
Let $K$ be a knot in $S^3$, and $E_K$ be the exterior manifold. 
Then gluing $S^1\times E_K$ to the boundary of $X^\prime$ produces a manifold $X_K$.
The manifold $X_K$ has the following properties.
\begin{enumerate}
\item $X_K$ is homeomorphic to $E(n)$. 
\item Let $A_K(t)=a_0 +\sum a_j(t^j+t^{-j})$ be the symmetrized Alexander polynomial of the knot $K$. 
For two knots $K_1$ and $K_2$, if $A_{K_1}(t)\not\equiv A_{K_2}(t)$, then $X_{K_1}$ is not diffeomorphic to $X_{K_2}$. 
\item Any exotic smooth structure obtained by a knot surgery can not be constructed by using log transforms. 
\end{enumerate}
To prove (2), Seiberg-Witten invariants are used.
In particular, it is known that 
$$
\SW_{X_K}(c_{\spin}) = a_0\cdot\SW_{E(n)}(c_{\spin}).
$$
Therefore, if both of $\SW_{E(n)}(c_{\spin})$ and $a_0$ are not divisible by $3$, then $\SW_{X_K}(c_{\spin})$ is also not divisible by $3$.
Note that there are many tori in $E(n)$ which give different homology classes, and further surgeries on these tori give more smooth structures.

Now, we introduce the following family of smooth structures on $E(n)$, and prove \corref{cor:ell-z3} and the assertion (2) of \thmref{thm:ellz3}.
\begin{Definition}\label{def:up}
Let $\UU_{E(n),p}$ be the set of smooth structures on $E(n)$ which consists of
\begin{itemize}
\item the standard smooth structure $E(n)$,
\item $E(n)_{k,l}$ for some odd and coprime $k,l$, and
\item smooth structures obtained by operating knot surgery constructions along tori in $E(n)_{k,l}$ by using knots $K$ which satisfy $a_0\not\equiv 0$ mod $p$.
\end{itemize}
\end{Definition}
\proof[Proof of \corref{cor:ell-z3}]
For each smooth structure in $\UU_{E(n),3}$, $\SW_X(c_{\spin})\not\equiv 0$ mod $3$ if $c_{n-2}\not\equiv 0$ mod $3$. 
Therefore \thmref{thm:spin-z3} proves the corollary.
\endproof
\proof[Proof of the assertion (2) of \thmref{thm:ellz3}]
This is clear by the facts that $c_{2}=2\not\equiv 0$ mod $3$, and the classes indicated as ``NS'' in \tabref{tab:ellz3} have positive $m_+$ and $m_-$.\endproof

\begin{Remark}\label{rem:en-z3}
By similar arguments as in \subsecref{subsec:loclin-z3} and \subsecref{subsec:const-z3}, we can obtain results similar to \thmref{thm:ellz3} for $\Z_3$-actions on homotopy $E(n)$ of larger $n$.
For example, classifications of $\Z_3$-actions on $E(8)$ and $E(10)$ are summarized in \tabref{table:e8-z3} and \tabref{table:e10-z3} at the end of the paper.
There are several remarks on tables.
\begin{enumerate}
\item Since $c_{n-2}\not\equiv 0$ mod $3$ for $n=8$ and  $10$, we can use \corref{cor:ell-z3} to judge the nonsmoothability of each class for smooth structures in $\UU_{E(n),3}$.  
Classes which admit no smooth action for such smooth structures are indicated as ``NS'' in tables. 
\item For the class no.$43$ in \tabref{table:e10-z3} indicated as ``No REP'',  we can not construct a $G$-invariant unimodular form by the method in this subsection. 
Therefore, we do not have a locally linear action of this class at present. 
Except this class, each of  classes in tables can be realized by a locally linear action by the method in this subsection.  
\item Note that, for given $b_+^G$, at most two classes have possibility to be smooth for above smooth structures.  
\end{enumerate}
\end{Remark}
\begin{Remark}
\tabref{table:en-z3} below summarizes numbers of classes of $\Z_3$-actions on $E(n)$.
In the table, the column ``Total'' is for the total numbers of classes for the classification by fixed point data.
The column ``NS'' is for the numbers of classes which turned out to be nonsmoothable for smooth structures in $\UU_{E(n),3}$ by our method using the mod $p$ vanishing theorem.
The column ``No REP'' is for the numbers of classes of which a locally linear action can not be constructed by our method. 
Note that, for $E(n)$, the ratio of ``NS'' classes in the total classes increases as $n$ increases.
In the case of $E(28)$, the ratio reaches $90$ percent.
\remref{rem:en-z3}, (3) can be considered as a reason.
\end{Remark}
\begin{table}[h]
\caption{Numbers of classes of $\Z_3$-actions on $E(n)$}\label{table:en-z3}
\begin{center}
\begin{tabular}{l|r|r|c}
& Total & NS & No REP\\
\hline
$E(2)$ & $4$ & $1$ & $0$\\
$E(4)$ & $10$ & $5$ & $1$\\
$E(8)$ & $30$ & $21$ & $0$\\
$E(10)$ & $44$ & $33$ & $1$\\
$E(20)$ & $154$ & $133$ & $0$\\
$E(22)$ & $184$ & $161$ & $1$\\
$E(26)$ & $252$ & $225$ & $0$\\
$E(28)$ & $290$ & $261$ & $1$\\
\end{tabular}
\end{center}
\end{table}

Finally, we prove the $\Z_3$-case of \thmref{thm:elln}.
\proof[Proof of the $\Z_3$-case of \thmref{thm:elln}]
For $X=E(n)$, we have $e=\chi(X)=12n$, $s=\Sign(X)=-8n$, $b_+=2n-1$ and $b_-=10n-1$. 
With these data, \eqref{eq:ineq-e}, \eqref{eq:mod3}, \eqref{eq:mod9-s} and \eqref{eq:ineq-bG} should be satisfied.
In fact, all possibilities of pairs $(m_+,m_-)$ can be written as
\begin{equation}\label{eq:Cn}
\CC(n):=
\left\{ (m_+,m_-)\in 3\Z\times 3\Z\,\left|
\begin{aligned}
k\in & \Z,\,  \frac13 n\leq k\leq n, \\
m_+&\geq 0, \, m_-\geq 0,\\ 
2m_+ &+m_-=9k-3n
\end{aligned}
\right.
\right\}.
\end{equation}
This $\CC(n)$ gives the classification table of $\Z_3$-actions on $E(n)$.
To obtain \eqref{eq:Cn}, consider as follows.
Since every nontrivial real representation of $G=\Z_3$ has even rank, $b_+^G$ can be written as $b_+^G=2k-1$, where $k$ is an integer which satisfies $1\leq k \leq n$.
>From \eqref{eq:bpG}, we obtain the relation $2m_+ +m_-=9k-3n$. 
Summing up all the other conditions with this, we can obtain \eqref{eq:Cn}.

For many (perhaps almost all) pairs in $\CC(n)$, we can construct corresponding $G$-invariant forms for $E(n)$ by the method in \subsecref{subsec:loclin-z3}.
Since it would be complicated to give a general procedure to construct $G$-forms for all pairs, we are content here to construct a $G$-form for a pair in $\CC(n)$ which gives a ``NS'' class.

Consider the pair $(m_+, m_-)=(\frac{3n}2, 3n)\in \CC(n)$. 
Then the corresponding $G$-form can be given as $(2n-1)A\oplus\frac{n}2 \Gamma_{16,5}$. 
Therefore, the conditions {\it REP} and {\it GSF} are satisfied, and we have a locally linear pseudofree $G$-action by \thmref{thm:EE}.

Since both $m_+$ and $m_-$ are positive, this $G$-action is nonsmoothable with respect to infinitely many smooth structures in $\UU_{E(n),3}$.
Thus the theorem is established.
\endproof
%
%
%
%
\section{$\Z_5$-actions on elliptic surfaces}\label{sec:z5}
%
%
In this section, we prove \thmref{thm:k3z5} and the $\Z_5$-case of \thmref{thm:elln}.
Since proofs are similar to those of the $\Z_3$-case in \secref{sec:z3}, details will be omitted.
%
%
\subsection{Existence of locally linear $\Z_5$-actions}\label{subsec:loclin-z5}
%
%
In this subsection, we prove the assertion (1) of \thmref{thm:k3z5}. 
The argument is  parallel to \subsecref{subsec:loclin-z3}.

There are six types of representations at fixed points for  $\Z_5$-actions. 
\begin{itemize}
\item The type ($11$): $(1,1)$ or $(4,4)$
\item The type ($22$): $(2,2)$ or $(3,3)$
\item The type ($12$): $(1,2)$ or $(3,4)$
\item The type ($13$): $(1,3)$ or $(2,4)$
\item The type ($14$): $(1,4)$ 
\item The type ($23$): $(2,3)$ 
\end{itemize}

Let $m_{ij}$ be the number of fixed points of the type ($ij$). 
Pseudofree locally linear $G$-actions have the following properties.
\begin{Proposition}\label{prop:propz5}
Suppose that $G=\Z_5$ acts on a simply-connected oriented closed manifold $X$. 
Let $e=\chi(X)$  and $s=\Sign (X)$. 
Then the data $\{m_{ij}\}$ satisfy the following{\rom :}
\begin{equation}\label{eq:ineq-e5}
m_{11} + m_{22} + m_{12} + m_{13} + m_{14} + m_{23}  \leq e.
\end{equation}
\begin{equation}\label{eq:mod5-e}
m_{11} + m_{22} + m_{12} + m_{13} + m_{14} + m_{23} \equiv -\frac{e}4  \mod 5.
\end{equation}
\begin{equation}\label{eq:mod5-s}
\left\{
\begin{aligned}
- m_{11} - m_{22} + m_{14} + m_{23} &\equiv -\frac{s}4 \mod 5,\\
- m_{11} + 3m_{22} - m_{12} + m_{13} + m_{14} - 3m_{23}&\equiv -s \mod 5,\\
3m_{11} - m_{22} + m_{12} - m_{13} - 3m_{14} + m_{23}&\equiv -s \mod 5.
\end{aligned}
\right.
\end{equation}
\begin{equation}\label{eq:ineq-bG5}
0\leq b_+^G \leq b_+,\, 0\leq b_-^G \leq b_-,
\end{equation}
where
\begin{align}\label{eq:z5bp}
b_+^G &= \frac15\left\{\frac{e+s}2+2m_{12} + 2m_{13} + 4m_{14} + 4m_{23})\right\}-1,\\\label{eq:z5bm}
b_-^G &= \frac15\left\{\frac{e-s}2+4m_{11} + 4m_{22}+ 2m_{12} + 2m_{13} )\right\}-1.
\end{align}
\end{Proposition}
\proof
By \thmref{thm:EE},the Lefschetz formula and the $G$-signature formula should be satisfied.
Calculations similar to \subsecref{subsec:loclin-z3} show the proposition.\endproof

By \propref{prop:propz5}, we can show that pseudofree locally linear $\Z_5$-actions on $K3$ are classified into $285$ classes of  the fixed point data $\{m_{ij}\}$. 
Let $\CC_{E(2),5}$ be the set of these $285$ classes.
\begin{Remark}
In our classification, {\it weakly equivalent} classes are identified: 
Suppose a $G=\Z_5$-action is given. 
One can consider another $G$-action which is given by $(g,x)\to g^2x$.
These two actions may have different $m_{ij}$, for instance, $m_{11}$ and $m_{22}$ exchange their values. 
However, these two are identified, since they are essentially same.
\end{Remark}

To prove the existence of a locally linear $G$-action for given data $\{m_{ij}\}$, we need to construct a $G$-form. 
\begin{Lemma}\label{lem:even5}
Let $r=16(5q+1)$, where $q$ is a non-negative integer.
For each integer $k$ which satisfies $0\leq k\leq 16q+3$, there is a $G$-action on $\Gamma_{r}$ such that
$$
\Gamma_{r}\cong (r-5k)\Z\oplus k\Z[G]\text{ as a $\Z[G]$-module}.
$$
This unimodular form with the $G$-action is denoted by $\Gamma^5_{r,k}$.
\end{Lemma}
\proof
When $k=0$, it suffices to take the trivial $G$-action. 
Hence we suppose $k\geq 1$. 

Note that the symmetric group of degree $r$ acts on $\Gamma_{r}$ as permutations of components.
For a fixed generator $g$ of $G$, define the $G$-action on $\Gamma_{r}$ by
$$g = (1,2,3,4,5)(6,7,8,9,10)\cdots(5k-4,5k-3,5k-2,5k-1,5k),$$
where $(l,m,n,o,p)$ is the cyclic permutation of $(x_l, x_m ,x_n,x_o,x_p)$.

As a basis for $\Gamma_{r}$, we take 
\begin{align*}
f_i &= \left\{
\begin{aligned}
e_i+&e_{r}, \quad\qquad (i=1,\ldots ,40q+10), \\
e_i-&e_{r}, \quad\qquad (i=40q+11,\ldots,r-2),
\end{aligned}\right.\\
f_{r-1} &=e_{r-1}-3e_{r}, \\
f_{r} &=\frac12 (e_1+e_2+\cdots+e_{r}),
\end{align*}
where $e_1,\ldots,e_{r}$ is the usual orthonormal basis for $\R^{r}$. 
Then the basis $(f_1,f_2,\ldots,f_{r})$ gives the required direct splitting if $k\leq 16q+2$.

When $k=16q+3$, we need to change basis.
New basis $\{f_i^\prime\}$ is given as follows:
$f^\prime_i = f_i$ for $i=1,\ldots ,r-6$, $f_r^\prime=f_r$ and
$(f_{r-5}^\prime,f_{r-4}^\prime,f_{r-3}^\prime,f_{r-2}^\prime,f_{r-1}^\prime)=(v,gv, g^2v,g^3v,g^4v)$ where $v=-f_{r-2}-f_{r-1}$. 
Then this basis gives the required property.
\endproof
\begin{Lemma}\label{lem:H5}
There is a $G$-invariant form $B^5_{1,1}$ such that $B^5_{1,1}\cong 3H$ as a form, and $B^5_{1,1}\cong \Z\oplus\Z[G]$ as a $\Z[G]$-module, and $r_+^G=r_-^G=1$, where $r_+^G$ {\rom (}resp. $r_-^G${\rom )} denote the rank of the $G$-fixed part of a maximal positive {\rom (}resp. negative{\rom )} definite subspace of $B^5_{1,1}\otimes\R$.
\end{Lemma}
\proof
Let us consider the basis of  $\Z\oplus\Z[G]$ of the form $\{f,e,ge,g^2e,g^3e,g^4e\}$. 
With respect to this basis, let us consider the form represented by the matrix
$$
Q=
\begin{pmatrix}
 2 & 1 & 1 & 1 & 1 & 1\\ 
 1 & 0 & 1 & 0 & 0 & 1\\
 1 & 1 & 0 & 1 & 0 & 0\\
 1 & 0 & 1 & 0 & 1 & 0\\
 1 & 0 & 0 & 1 & 0 & 1 \\
 1 & 1 & 0 & 0 & 1 & 0
\end{pmatrix}.
$$
It is easy to see that $Q$ represents a $G$-invariant unimodular even form on $\Z\oplus\Z[G]$ and that $r_+^G=r_-^G=1$. 
Since the rank of $Q$ is $6$, $Q$ should be equivalent to $3H$ as a form. 
\endproof

Let $A^5$ be the $G$-invariant unimodular form such that $A^5\cong H$ as a form and $A^5\cong 2\Z$ as a $\Z[G]$-module.

To each class in $\CC_{E(2),5}$, we associate the following form:
\begin{itemize}
\item $3A^5\oplus\Gamma^5_{16,k}$ to classes with $b_+^G=3$ and $b_-^G=19-4k$,
\item $B^5_{1,1}\oplus\Gamma^5_{16,k}$ to classes with $b_+^G=1$ and $b_-^G=17-4k$.
\end{itemize}
With these forms, we can prove that every class except one class in $\CC_{E(2),5}$ can be realized by a locally linear actions by \thmref{thm:EE}. 
The only one exception is the class given by $m_{14}=m_{23}=2$ and $m_{11}=m_{22}=m_{12}=m_{13}=0$.
However, we can construct a smooth action of this class as follows.
\begin{Proposition}\label{prop:z5-1}
There exists a smooth action on the projective $K3$ surface in $\CP^4$ which satisfies $m_{14}=m_{23}=2$ and $m_{11}=m_{22}=m_{12}=m_{13}=0$.
\end{Proposition}
\proof
Let us consider the $K3$ surface $X$ defined by equations $\sum_{i=0}^4 z_i^2 =0$ and $\sum_{i=0}^4 z_i^3 =0$ in $\CP^4$. 
By the symmetry of defining equations, the symmetric group of degree $5$ acts on $X$ by permutations of variables.
Via this action, $G$ acts on $X$ smoothly (in fact, holomorphically). 
It is easy to see that this $G$-action has the required property.
\endproof
\begin{Remark}\label{rem:z5-502}
The class given by $m_{11}=1$, $m_{12}=3$ and $m_{22}=m_{13}=m_{14}=m_{23}=0$ is also realized by a smooth action.  
In fact, the holomorphic action in Example 5.4 of \cite{Zh} belongs to this class.
\end{Remark}
%
\subsection{A constraint on smooth $\Z_5$-actions on elliptic surfaces}\label{subsec:const-z5}
%
%
In this subsection, we prove a proposition which gives a constraint on smooth $\Z_5$-actions on elliptic surfaces, and finally prove the assertion (2) of \thmref{thm:k3z5} and the $\Z_5$-case of \thmref{thm:elln}.
%
%
%
%
\begin{Proposition}\label{prop:spin-z5}
Let $G=\Z_5$, and $X$ be a simply-connected closed oriented smooth spin $4$-manifold with $b_+\geq 2$, which satisfies $2\chi(X) + 3\Sign(X)=0$.
Suppose $G$ acts on $X$ smoothly and pseudofreely so that  $b_+^G\geq 1$. 
If $\SW_X(c_{\spin})\not\equiv 0$ mod $5$, 
then at least one of the following holds,
\begin{align}\label{eq:z5m1}
m_{11}&=m_{22}=m_{12}=m_{13}=0,\\\label{eq:z5m2}
\text{or }\quad
m_{22}&\geq 2m_{13}+2m_{14}+3m_{23},\\\label{eq:z5m3}
\text{or }\quad
m_{11}&\geq 2m_{12}+3m_{14}+2m_{23}.
\end{align}
\end{Proposition}
\proof
Let $G$ act on $c_{\spin}$ as in \remref{rmk:spin-action}, and take the trivial flat connection $A_0$ on $L$ as the reference $G$-invariant connection.
Then $\SW_X(c_{\spin})\not\equiv 0$ mod $5$ implies that there exist $j$ which satisfies $2k_j\geq 1+b_+^G$.

By the $G$-spin theorem, coefficients $k_j$ are given as follows:
\begin{align*}
k_0 &=\frac15 \left(-\frac{s}8-2m_{11}-2m_{22}+2m_{14}+2m_{23}\right),\\
k_1 = k_4 &= \frac15 \left(-\frac{s}8 +m_{22}+m_{12}-m_{13}-m_{23}\right),\\
k_2 = k_3 &= \frac15 \left(-\frac{s}8 +m_{11}-m_{12}+m_{13}-m_{14}\right).
\end{align*}
Note that $b_+^G$ has already calculated in \eqref{eq:z5bp}. 
By using the relation $2e+3s=0$, we can show that $2k_0\geq 1+b_+^G$, $2k_1=2k_4\geq 1+b_+^G$ and $2k_2=2k_3\geq 1+b_+^G$ are equivalent to \eqref{eq:z5m1}, \eqref{eq:z5m2} and \eqref{eq:z5m3}, respectively.
\endproof
\begin{Corollary}\label{cor:en5}
Let $X$ be a homotopy $E(n)$ of even and positive $n$ with a smooth structure in $\UU_{E(n),5}$.
Suppose $G=\Z_5$ act on $X$ smoothly and pseudofreely so that $b_+^G\geq1$.
If $c_{n-2}\not\equiv 0$ mod $5$, then at least one of \eqref{eq:z5m1}, \eqref{eq:z5m2} or \eqref{eq:z5m3} holds.
\end{Corollary}
\proof[Proof of the assertion (2) of \thmref{thm:k3z5}]
This is obvious because $c_0=1$ and \corref{cor:en5}.\endproof
Now, we prove the $\Z_5$-case of \thmref{thm:elln}.
\proof[Proof of the $\Z_5$-case of \thmref{thm:elln}]
We can construct many examples of  classes for each $n$ which do not satisfy any of \eqref{eq:z5m1}, \eqref{eq:z5m2} and \eqref{eq:z5m3}.
Since it would be complicated to give a general procedure to enumerate all of such classes, we are content here to construct such an example for each $n$.

First note that $c_{n-2}\equiv 0$ mod $5$ if $n\equiv 0$ or $8$ mod $10$. 
(See \remref{rem:p-adic}.)
Hence we assume $n\not\equiv 0$ or $8$ mod $10$.
Then, examples as above are given according to $n$ as follows:
Let $l$ be a non-negative integer.
\begin{enumerate}
\item When $n=10l+2$, $m_{22}=1$, $m_{13}=40l+3$, $m_{12}=5$, $m_{11}=m_{14}=m_{23}=0$.
\item When $n=10l+4$, $m_{22}=2$, $m_{13}=40l+11$, $m_{12}=5$, $m_{11}=m_{14}=m_{23}=0$.
\item When $n=10l+6$, $m_{22}=3$, $m_{13}=40l+19$, $m_{12}=5$, $m_{11}=m_{14}=m_{23}=0$.
\end{enumerate}
It is easy to check that these do not satisfy any of \eqref{eq:z5m1}, \eqref{eq:z5m2} and \eqref{eq:z5m3}.

These satisfy conditions \eqref{eq:ineq-e5}--\eqref{eq:ineq-bG5}. 
Furthermore, we can construct a $G$-form corresponding to each class above as follows.
\begin{enumerate}
\item $(20l+3)A^5 \oplus \Gamma_{16(5l+1),16l+3}^5.$
\item $(20l+7)A^5 \oplus \Gamma_{16(5l+1),16l+3}^5\oplus\Gamma^5_{16,3}.$
\item $(20l+11)A^5 \oplus \Gamma_{16(5l+1),16l+3}^5\oplus\Gamma^5_{16,3}\oplus\Gamma^5_{16,3}.$
\end{enumerate}
(For forms $A^5$ and $\Gamma_{r,k}^5$, see \lemref{lem:even5} and \lemref{lem:H5} below.)
Therefore, the conditions {\it REP} and {\it GSF} are satisfied, and we have a locally linear pseudofree $G$-action for each class above by \thmref{thm:EE}.

Then, \corref{cor:en5} implies that each of these $G$-actions is nonsmoothable with respect to infinitely many smooth structures in $\UU_{E(n),5}$.
Thus the theorem is established.
\endproof
%
%
%
%
\section{$\Z_7$-actions on elliptic surfaces}\label{sec:z7}
%
%
In this section, we prove the $\Z_7$-case of \thmref{thm:elln}. 
%
%
\subsection{Existence of locally linear $\Z_7$-actions}\label{subsec:loclin-z7}

There are twelve types of representations at fixed points for $G=\Z_7$-actions: 
$(1,1)$, $(2,2)$, $(3,3)$, 
$(1,2)$, $(2,4)$, $(1,4)$, 
$(1,5)$, $(2,3)$, $(1,3)$, 
$(1,6)$, $(2,5)$, $(3,4)$.
Let $m_{ij}$ be the number of fixed points of the type $(i,j)$.

\begin{Proposition}
Suppose that $G$ acts on a simply-connected closed oriented $4$-manifold $X$ locally linearly and pseudofreely.
Let $e=\chi(X)$ and $s=\Sign(X)$.
Then the data $\{m_{ij}\}$ satisfy
\begin{equation}\label{eq:z7-sum}
\sum m_{ij} \leq e,\quad \sum m_{ij} \equiv -6e \mod 7,
\end{equation}
\begin{multline}
-10 (m_{11} + m_{22} + m_{33})  -2 ( m_{12} + m_{24} + m_{14}) \\
+10 (m_{16} + m_{25} +m_{34}  ) 
+2( m_{15} + m_{23} + m_{13})\equiv -s  \mod 7,
\end{multline}
\begin{multline}
-5  m_{11} + 7 m_{22} + 3 m_{33} - 3 m_{12} + 3 m_{24} + m_{14} \\
+5  m_{16} - 7 m_{25} - 3 m_{34} + 3 m_{15} - 3 m_{23} - m_{13}\equiv -s \mod 7,
\end{multline}
\begin{multline}
 3  m_{11} - 5 m_{22} + 7 m_{33} +   m_{12} - 3 m_{24} +3m_{14}\\
- 3  m_{16} + 5 m_{25} - 7 m_{34} -   m_{15} + 3 m_{23} -3m_{13}\equiv -s \mod 7,
\end{multline}
\begin{multline}\label{eq:l7-4}
 7  m_{11} + 3 m_{22} - 5 m_{33} + 3 m_{12} +   m_{24} -3m_{14}\\
- 7  m_{16} - 3 m_{25} + 5 m_{34} - 3 m_{15} -   m_{23} +3m_{13}\equiv -s \mod 7,
\end{multline}
\begin{equation}\label{eq:z7-b}
0\leq b_+^G\leq b_+,\quad 0\leq b_-^G\leq b_-,
\end{equation}
where
\begin{equation}
\begin{split}
b_+^G =\frac17 \left\{ \frac{e + s}2  - 2 \right.&(m_{11} + m_{22} + m_{33})  +2 ( m_{12} + m_{24} + m_{14}) \\ 
+ 8 &(m_{16} + m_{25} +m_{34}  ) +4 ( m_{15} + m_{23} + m_{13})\bigg\}-1,
\end{split}
\end{equation}
\begin{equation}
\begin{split}
b_-^G =\frac17 \left\{ \frac{e - s}2   \right.
+ 8 &(m_{11} + m_{22} + m_{33})  +4 ( m_{12} + m_{24} + m_{14}) \\ 
- 2 &(m_{16} + m_{25} +m_{34}  ) +2 ( m_{15} + m_{23} + m_{13})\bigg\}-1.
\end{split}
\end{equation}
\end{Proposition} 
\proof
By \thmref{thm:EE}, the Lefschetz formula and the $G$-signature formula should be satisfied. 
Calculations similar to \subsecref{subsec:loclin-z3} show the proposition. 
\endproof
Conversely, the data $\{m_{ij}\}$ which satisfies \eqref{eq:z7-sum}--\eqref{eq:z7-b} can be realized as fixed point data of a locally linear action if a corresponding $G$-form is constructed.
In particular, we can always obtain locally linear actions in the homologically trivial cases:
\begin{Proposition}\label{prop:z7-ht}
Let $X$ be $E(n)$ of even positive $n$.
If $\{m_{ij}\}$ which satisfies \eqref{eq:z7-sum}--\eqref{eq:l7-4} and $b_+^G=b_+$ and $b_-^G=b_-$, then the data $\{m_{ij}\}$  can be realized as fixed point data of a locally linear action.
\end{Proposition}
\proof 
Consider the trivial $G$-form. Then  \thmref{thm:EE} proves the proposition.
\endproof
\begin{Remark}
In order to construct other types of $G$-forms, we can prove a lemma similar to \lemref{lem:even} or \lemref{lem:even5}. 
However, in the case of $\Z_7$-actions, we have a plenty of classes of homologically trivial actions, and it will turn out to suffice to consider such classes for our purpose. 
Therefore, we do not write down the lemma for such $G$-forms here. 
\end{Remark}
%
%
%
%
\subsection{A constraint on smooth $\Z_7$-actions on elliptic surfaces}\label{subsec:const-z7}
%
%
This subsection proves a proposition which gives a constraint on smooth $\Z_7$-actions on elliptic surfaces, and finally proves the $\Z_7$-case of \thmref{thm:elln}.
\begin{Proposition}\label{prop:spin-z7}
Let $G=\Z_7$, and $X$ be a simply-connected closed oriented smooth spin $4$-manifold with $b_+\geq 2$ which satisfies $2\chi(X) + 3\Sign(X)=0$. 
Suppose $G$ acts on $X$ smoothly and pseudofreely so that $b_+^G\geq 1$.
If $\SW_X(c_{\spin})\not\equiv 0$ mod $7$, then at least one of the following holds,
\begin{align}\label{eq:z7m1}
m_{12}&+m_{24}+m_{14}\geq 3(m_{11}+m_{22}+m_{33})+4(m_{15}+m_{23}+m_{13}),\\\label{eq:z7m2}
\text{or }\quad
3m_{22}&+2m_{33}\geq m_{12}+3m_{24}+3m_{16}+6m_{25}+5m_{34}+2m_{15}+3m_{13},\\\label{eq:z7m3}
\text{or }\quad
2m_{11}&+3m_{33}\geq m_{24}+3m_{13}+5m_{16}+3m_{25}+6m_{34}+3m_{15}+2m_{23},\\\label{eq:z7m4}
\text{or }\quad
3m_{11}&+2m_{22}\geq 3m_{12}+ m_{14}+6m_{16}+5m_{25}+3m_{34}+3m_{23}+2m_{13}.
\end{align}
\end{Proposition}
\proof
The proof is similar to that of \propref{prop:spin-z5}. 
The coefficients $k_j$ of the $G$-index of the Dirac operator are calculated from the $G$-spin formula as follows.
\begin{equation}\label{eq:z7k0}
\begin{split}
k_0 = \frac17 \left\{ -\frac18 s\right.   
-4&(m_{11} + m_{22} + m_{33})  +2 ( m_{12} + m_{24} + m_{14})\\ 
+4&(m_{16} + m_{25} +m_{34}  ) -2 ( m_{15} + m_{23} + m_{13})\bigg\},
\end{split}
\end{equation}

\begin{equation}\label{eq:z7k1}
\begin{split}
k_1 = k_6 = \frac17 \left\{ -\frac18 s\right.    
-  & m_{11} + 2 m_{22} +   m_{33}  - 2 m_{24} + m_{14}\\
+  & m_{16} - 2 m_{25} -   m_{34}  + 2 m_{23} - m_{13}\bigg\},
\end{split}
\end{equation}

\begin{equation}\label{eq:z7k2}
\begin{split}
k_2 = k_5 = \frac17 \left\{-\frac18 s\right.    
+   & m_{11} -   m_{22} + 2 m_{33} +   m_{12}  -2m_{14}\\
-   & m_{16} +   m_{25} - 2 m_{34} -   m_{15}  +2m_{13}\bigg\},
\end{split}
\end{equation}

\begin{equation}\label{eq:z7k3}
\begin{split}
k_3 = k_4 = \frac17 \left\{-\frac18 s\right.    
+ 2 & m_{11} +   m_{22} -   m_{33} - 2 m_{12} +   m_{24} \\
- 2 & m_{16} -   m_{25} +   m_{34} + 2 m_{15} -   m_{23} \bigg\},
\end{split}
\end{equation}
\thmref{thm:mod-p} implies that there exists $k_j$ so that $2k_j\geq 1+b_+^G$. 
The proposition is obtained by rewriting these inequalities in terms of $m_{ij}$.
\endproof
Now, let us prove the $\Z_7$-case of \thmref{thm:elln}.
\proof[Proof of the $\Z_7$-case of \thmref{thm:elln}]
Consider the data $m_{33}=m_{24}=7n/2$, $m_{14}=4n$, $m_{23}=n$ and all other $m_{ij}$ are zero. 
Since these $\{m_{ij}\}$ satisfy \eqref{eq:z7-sum}--\eqref{eq:l7-4}, $b_+^G=b_+$ and $b_-^G=b_-$, there exists a locally linear $G$-action with fixed point data $\{m_{ij}\}$ by \propref{prop:z7-ht}. 
On the other hand, these $\{m_{ij}\}$ do not satisfy any of \eqref{eq:z7m1}--\eqref{eq:z7m4}. 
Therefore the $G$-action is nonsmoothable with respect to the infinitely many smooth structures in $\UU_{E(n),7}$.
Thus, \thmref{thm:elln} is established.
\endproof
\begin{Remark}
In the case of $K3$, there are $124256$ classes of $\Z_7$-actions. 
(In this enumeration, weakly equivalent classes are identified.)
Among these, $103829$ classes are homologically trivial. 
Therefore, these classes can be realized by locally linear actions by \propref{prop:z7-ht}.
On the other hand, there are $4772$ homologically trivial classes do not satisfy \eqref{eq:z7m1}--\eqref{eq:z7m4}. 
Therefore, by \propref{prop:spin-z7}, locally linear actions in these classes are nonsmoothable for smooth structures in $\UU_{E(2),7}$. 
\end{Remark}
\begin{Remark}
The authors know only one example of smooth $\Z_7$-action on $K3$ which is given in \cite{Zh}(Example 5.4).
The fixed point data of this action is as follows: $m_{16}=1$, $m_{13}=2$ and all other $m_{ij}$ are zero.
This action is not homologically trivial.
\end{Remark}
%
%
\section{Concluding remarks}\label{sec:concluding}
%
%
In this last section, we give several remarks.
%
%
\subsection{Other $\Spinc$-structures}
%
%
In the arguments so far, we use only $c_{\spin}$, while ellptic surfaces except the standard $K3$ have basic classes other than $c_{\spin}$.
By using such basic classes, we can obtain more constraints on smooth actions.

For example, let $X$ be the standard $E(4)$, and $G=\Z_3$. 
According to \thmref{thm:ellz3}, a smooth $G$-action on $X$ may have the data $m_+=12 $ and $m_-=0$.
\begin{Theorem}\label{thm:pdf}
Suppose $G=\Z_3$ acts on $X=E(4)$ smoothly and pseudofreely so that $m_+=12 $ and $m_-=0$.
Let $PD[F]$ be the Poincar\'{e} dual of the homology class of a regular fiber $F$, $L$ the complex line bundle whose $c_1$ is $PD[F]$, and $\pi\colon X\to X/G$ the projection to the quotient space.
Then, there exists a complex line bundle $\bar{L}$ over the quotient space $X/G$ such that $L=\pi^*\bar{L}$.
\end{Theorem}
To prove \thmref{thm:pdf}, first note the following lemma.
\begin{Lemma}Under the assumption of \thmref{thm:pdf}, the $G$-action lifts to the $\Spinc$-structure $c$ such that $\det c =L$.
\end{Lemma}
\proof
Note that, if $\SW_X(c)=1$, then $\det c =\pm PD[F]$.
Since the $G$-action preserves basic classes, $g(PD[F])=\pm PD[F]$ for any $g\in G$.
Since $g^3 =1$, $g(PD[F])=PD[F]$ for any $g\in G$.
Thus, the $G$-action preserves $c_1(L)=PD[F]$, and by the theorem of Hattori-Yoshida\cite{HY}, the $G$-action lifts to $L$. 
Then, the $G$-action lifts to $c$, since $G$ is odd order.
\endproof
Thus the $G$-action lifts to $L$ and $c$, and $G$ acts on the fiber of $L$ over each fixed point with some weight.
Then the mod $p$ vanishing theorem implies
\begin{Lemma}
Every fixed point has the same weight on $L$.
\end{Lemma}
Therefore, $L$ can be considered as the pull-back of a line bundle $\bar{L}$ on $X/G$. 
Thus, \thmref{thm:pdf} is proved.
Note that we can obtain similar results for other situations, i.e., other fixed point data, other $G$ and other $E(n)$.
%
%
\subsection{Dependence on smooth structures}\label{subsec:dep}
%
%
As mentioned in the introduction, many authors have constructed a lot of examples of nonsmoothable locally linear actions \cite{KLee,KLaw, HL, Bryan, Kiyono, ChenKw}. 
In the papers \cite{KLaw, HL,Bryan,Kiyono, ChenKw}, the authors use gauge theory to prove that the actions are nonsmoothable. 
It is interesting that the actions in \cite{KLaw,HL,Bryan,Kiyono} are nonsmoothable for {\it arbitrary} smooth structures:
In \cite{Bryan} and \cite{Kiyono}, Bryan and Kiyono use some $G$-equivariant variants of $10/8$-inequalities which give constraints on $b_2$ and signature which do not depend on smooth structures.
Therefore the locally linear actions which violate these inequalities are clearly nonsmoothable  for {\it arbitrary} smooth structures.

On the other hand, in our case, we need to check the Seiberg-Witten invariant for {\it each} smooth structure in order to judge the nonsmoothability. 
This fact would suggest that our examples could be {\it subtle} in that the smoothablity of each action might depend on smooth structures.

At present, such subtle examples are known only by Chen-Kwasik \cite{ChenKw}. 
In \cite{ChenKw}, Chen and Kwasik prove that there is a family of symplectic exotic $K3$ surfaces  on which every nontrivial $\Z_p$-action of prime $p\geq 7$  is nonsmoothable, while there exist several examples of smooth $\Z_p$-actions of such $p$ on the standard $K3$.
This means that there are locally linear actions on $K3$ whose smoothabilities depend on smooth structures.

With these understood, the following problem would be interesting. 
({\it cf.} \thmref{thm:spin-z3} and \corref{cor:ell-z3}).
\begin{Problem}
Let $n$ be an even positive integer such that $c_{n-2}$ in \eqref{eq:cn} satisfies $c_{n-2}\not\equiv 0$ mod $3$.
Is there a smooth structure on $E(n)$ which admits a smooth $\Z_3$-action with $m_+>0$ and $m_->0$?
\end{Problem} 
To attack this problem, one could try to construct such a smooth $\Z_3$-action on a manifold  $X_K$ obtained from a knot surgery by a knot $K$ with $a_0\equiv 0$ mod $3$. 
(In this case, $\SW_{X_K}(c_{\spin})\equiv 0$ mod $3$. )

As the final remark, we note that, although we stop our calculations up to $p\leq 7$, nonsmoothable actions of higher order cyclic groups would be found by the same method.

\tiny
    \begin{table}[p]
\begin{minipage}{8.0cm}
    \caption{$\Z/3$-actions on $E(8)$}\label{table:e8-z3}
    \begin{center}
    \begin{tabular}{l|c|c|c|c|c|c|c}
    Class & $\#X^G$ & $m_+$ & $m_-$ & $b_2^G$ & $b_+^G$ & $b_-^G$ &\\
\hline
\hline
  1 & 24 & 24 &  0 & 46 & 15 & 31 &\\
  2 & 27 & 21 &  6 & 48 & 15 & 33 &NS\\
  3 & 30 & 18 & 12 & 50 & 15 & 35 &NS\\
  4 & 33 & 15 & 18 & 52 & 15 & 37 &NS\\
  5 & 36 & 12 & 24 & 54 & 15 & 39 &NS\\
  6 & 39 &  9 & 30 & 56 & 15 & 41 &NS\\
  7 & 42 &  6 & 36 & 58 & 15 & 43 &NS\\
  8 & 45 &  3 & 42 & 60 & 15 & 45 &NS\\
  9 & 48 &  0 & 48 & 62 & 15 & 47 &\\
\hline
 10 & 21 & 18 &  3 & 44 & 13 & 31 &NS\\
 11 & 24 & 15 &  9 & 46 & 13 & 33 &NS\\
 12 & 27 & 12 & 15 & 48 & 13 & 35 &NS\\
 13 & 30 &  9 & 21 & 50 & 13 & 37 &NS\\
 14 & 33 &  6 & 27 & 52 & 13 & 39 &NS\\
 15 & 36 &  3 & 33 & 54 & 13 & 41 &NS\\
 16 & 39 &  0 & 39 & 56 & 13 & 43 &\\
\hline
 17 & 15 & 15 &  0 & 40 & 11 & 29 &\\
 18 & 18 & 12 &  6 & 42 & 11 & 31 &NS\\
 19 & 21 &  9 & 12 & 44 & 11 & 33 &NS\\
 20 & 24 &  6 & 18 & 46 & 11 & 35 &NS\\
 21 & 27 &  3 & 24 & 48 & 11 & 37 &NS\\
 22 & 30 &  0 & 30 & 50 & 11 & 39 &\\
\hline
 23 & 12 &  9 &  3 & 38 &  9 & 29 &NS\\
 24 & 15 &  6 &  9 & 40 &  9 & 31 &NS\\
 25 & 18 &  3 & 15 & 42 &  9 & 33 &NS\\
 26 & 21 &  0 & 21 & 44 &  9 & 35 &\\
\hline
 27 &  6 &  6 &  0 & 34 &  7 & 27 &\\
 28 &  9 &  3 &  6 & 36 &  7 & 29 &NS\\
 29 & 12 &  0 & 12 & 38 &  7 & 31 &\\
\hline
 30 &  3 &  0 &  3 & 32 &  5 & 27 &\\
\end{tabular}
\end{center}
\end{minipage}
\begin{minipage}{7.1cm}
    \caption{$\Z/3$-actions on $E(10)$}\label{table:e10-z3}
    \begin{center}
    \begin{tabular}{l|c|c|c|c|c|c|c}
    Class & $\#X^G$ & $m_+$ & $m_-$ & $b_2^G$ & $b_+^G$ & $b_-^G$ &\\
\hline
\hline
  1 & 30 & 30 &  0 & 58 & 19 & 39 &\\
  2 & 33 & 27 &  6 & 60 & 19 & 41 &NS\\
  3 & 36 & 24 & 12 & 62 & 19 & 43 &NS\\
  4 & 39 & 21 & 18 & 64 & 19 & 45 &NS\\
  5 & 42 & 18 & 24 & 66 & 19 & 47 &NS\\
  6 & 45 & 15 & 30 & 68 & 19 & 49 &NS\\
  7 & 48 & 12 & 36 & 70 & 19 & 51 &NS\\
  8 & 51 &  9 & 42 & 72 & 19 & 53 &NS\\
  9 & 54 &  6 & 48 & 74 & 19 & 55 &NS\\
 10 & 57 &  3 & 54 & 76 & 19 & 57 &NS\\
 11 & 60 &  0 & 60 & 78 & 19 & 59 &\\
\hline
 12 & 27 & 24 &  3 & 56 & 17 & 39 &NS\\
 13 & 30 & 21 &  9 & 58 & 17 & 41 &NS\\
 14 & 33 & 18 & 15 & 60 & 17 & 43 &NS\\
 15 & 36 & 15 & 21 & 62 & 17 & 45 &NS\\
 16 & 39 & 12 & 27 & 64 & 17 & 47 &NS\\
 17 & 42 &  9 & 33 & 66 & 17 & 49 &NS\\
 18 & 45 &  6 & 39 & 68 & 17 & 51 &NS\\
 19 & 48 &  3 & 45 & 70 & 17 & 53 &NS\\
 20 & 51 &  0 & 51 & 72 & 17 & 55 &\\
\hline
 21 & 21 & 21 &  0 & 52 & 15 & 37 &\\
 22 & 24 & 18 &  6 & 54 & 15 & 39 &NS\\
 23 & 27 & 15 & 12 & 56 & 15 & 41 &NS\\
 24 & 30 & 12 & 18 & 58 & 15 & 43 &NS\\
 25 & 33 &  9 & 24 & 60 & 15 & 45 &NS\\
 26 & 36 &  6 & 30 & 62 & 15 & 47 &NS\\
 27 & 39 &  3 & 36 & 64 & 15 & 49 &NS\\
 28 & 42 &  0 & 42 & 66 & 15 & 51 &\\
\hline
 29 & 18 & 15 &  3 & 50 & 13 & 37 &NS\\
 30 & 21 & 12 &  9 & 52 & 13 & 39 &NS\\
 31 & 24 &  9 & 15 & 54 & 13 & 41 &NS\\
 32 & 27 &  6 & 21 & 56 & 13 & 43 &NS\\
 33 & 30 &  3 & 27 & 58 & 13 & 45 &NS\\
 34 & 33 &  0 & 33 & 60 & 13 & 47 &\\
\hline
 35 & 12 & 12 &  0 & 46 & 11 & 35 &\\
 36 & 15 &  9 &  6 & 48 & 11 & 37 &NS\\
 37 & 18 &  6 & 12 & 50 & 11 & 39 &NS\\
 38 & 21 &  3 & 18 & 52 & 11 & 41 &NS\\
 39 & 24 &  0 & 24 & 54 & 11 & 43 &\\
\hline
 40 &  9 &  6 &  3 & 44 &  9 & 35 &NS\\
 41 & 12 &  3 &  9 & 46 &  9 & 37 &NS\\
 42 & 15 &  0 & 15 & 48 &  9 & 39 &\\
\hline
 43 &  3 &  3 &  0 & 40 &  7 & 33 &No REP\\
 44 &  6 &  0 &  6 & 42 &  7 & 35 &\\
\end{tabular}
\end{center}
\end{minipage}
\end{table}

\end{document}